\documentclass[11pt]{article}
\usepackage{mathrsfs}
\usepackage{amssymb}
\usepackage{amsmath}
\usepackage{latexsym}
\usepackage{amsfonts}

\usepackage[all]{xy}
\def\Ker{\mathop{\rm Ker}\nolimits}
\def\Ext{\mathop{\rm Ext}\nolimits}
\def\Hom{\mathop{\rm Hom}\nolimits}
\def\End{\mathop{\rm End}\nolimits}
\def\Coker{\mathop{\rm Coker}\nolimits}
\def\Ind{\mathop{\rm Ind}\nolimits}
\def\mod{\mathop{\rm mod}\nolimits}

\def\id{\mathop{\rm id}\nolimits}
\def\pd{\mathop{\rm pd}\nolimits}
\def\gldim{\mathop{\rm gl.dim}\nolimits}
\def\add{\mathop{\rm add}\nolimits}

\title{\Large \bf Trivial Maximal 1-Orthogonal Subcategories For\\
Auslander's 1-Gorenstein Algebras
}
\author{Zhaoyong Huang\thanks{{\it E-mail address}:
huangzy@nju.edu.cn}, Xiaojin Zhang\thanks{{\it E-mail address}:
xiaojinzhang@sohu.com}\\
{\footnotesize \it Department of Mathematics, Nanjing University,
Nanjing 210093, Jiangsu Province, P. R. China}\\\\
{\small Dedicated to Professor Boxun Zhou on the occasion of his
ninetieth birthday}}
\date{ }
\begin{document}
\baselineskip=18pt
\maketitle

\vspace{0.2cm}

\begin{abstract}
Let $\Lambda$ be an Auslander's 1-Gorenstein Artinian algebra with
global dimension two. If $\Lambda$ admits a trivial maximal
1-orthogonal subcategory of $\mod\Lambda$, then for any
indecomposable module $M \in \mod \Lambda$, we have that the
projective dimension of $M$ is equal to one if and only if so is its
injective dimension and that $M$ is injective if the projective
dimension of $M$ is equal to two. In this case, we further get that
$\Lambda$ is a tilted algebra.
\end{abstract}

\vspace{0.5cm}

\centerline{\large \bf 1. Introduction}

\vspace{0.2cm}

Iyama in [Iy1] introduced the notion of maximal $n$-orthogonal
subcategories for Artinian algebras, and developed an extensive
theory. Maximal $n$-orthogonal subcategories are also called
$(n+1)$-cluster tilting subcategories following Keller and Reiten's
terminology in [KR]. It is an interesting problem to investigate
when maximal $n$-orthogonal subcategories (or modules) exist,
because it is the setting where Iyama has developed his higher
theory of almost split sequences and Auslander algebras (see [Iy1]
and [Iy2] for the details). The problem of the existence of maximal
$n$-orthogonal subcategories has been studied by several authors in
[EH], [GLS], [HZ1], [HZ2], [Iy1], [Iy2], [Iy3], and so on.

Iyama proved in [Iy3] that if $\Lambda$ is a finite-dimensional
algebra of finite representation type with Auslander algebra
$\Gamma$ and $\mod \Gamma$ contains a maximal 1-orthogonal module,
then $\Lambda$ is hereditary and the injective envelope of
$_{\Lambda}\Lambda$ is projective; furthermore, if the base field is
algebraically closed, then such $\Lambda$ is an upper triangular
matrix ring. On the other hand, we gave in [HZ1, Corollary 3.10] a
necessary and sufficient condition for judging when an Auslander
algebra admits a trivial maximal 1-orthogonal subcategory in terms
of the homological dimensions of simple modules as follows. For an
Auslander algebra $\Lambda$ with global dimension two, $\Lambda$
admits a trivial maximal 1-orthogonal subcategory of $\mod\Lambda$
if and only if any simple module $S \in \mod \Lambda$ with
projective dimension two is injective. This result played an
important role in studying the properties of Auslander algebras
admitting a trivial maximal 1-orthogonal subcategory. For example,
by applying it, we showed in [HZ2, Corollary 4.12] the following
result. Let $\Lambda$ be an Auslander algebra with global dimension
two. If $\Lambda$ admits a trivial maximal 1-orthogonal subcategory
of $\mod\Lambda$, then $\Lambda$ is a tilted algebra of finite
representation type. We also gave in [HZ2] an example to illustrate
that the conditions ``$\Lambda$ is an Auslander algebra" and
``$\Lambda$ is a tilted algebra of finite representation type" in
this result cannot be exchanged. In this paper, we study the
existence of trivial maximal 1-orthogonal subcategories under a
weaker condition than above, and prove the following result.

\vspace{0.2cm}

{\bf Theorem} {\it Let $\Lambda$ be an Auslander's 1-Gorenstein
Artinian algebra with global dimension two. If $\Lambda$ admits a
trivial maximal 1-orthogonal subcategory of $\mod\Lambda$, then for
any indecomposable module $M\in\mod\Lambda$, we have

(1) The projective dimension of $M$ is equal to one if and only if
so is its injective dimension.

(2) $M$ is injective if the projective dimension of $M$ is equal to
two.

Furthermore, we have that $\Lambda$ is a tilted algebra.}

\vspace{0.2cm}

Recall that an Artinian algebra $\Lambda$ is called {\it Auslander's
1-Gorenstein} if the injective envelope of $_{\Lambda}\Lambda$ is
projective. It was proved in [FGR, Theorem 3.7] that the notion of
Auslander's 1-Gorenstein algebras is left and right symmetric, that
is, $\Lambda$ is Auslander's 1-Gorenstein if and only if so is
$\Lambda ^{op}$. Also recall that an Artinian algebra $\Lambda$ is
called an {\it Auslander algebra} if the global dimension of
$\Lambda$ is at most two and the first two terms in a minimal
injective resolution of $_{\Lambda}\Lambda$ are projective. So it is
trivial that an Auslander algebra is Auslander's 1-Gorenstein.

This paper is organized as follows. In Section 2, we give some
notions and notations and collect some preliminary results. In
Section 3, we study the homological properties of simple modules and
indecomposable modules over an Artinian algebra with global
dimension two and admitting a trivial maximal 1-orthogonal
subcategory. As an application of the obtained properties, we prove
the above theorem. At the end of this section, we give an example to
illustrate this theorem.

\newpage

\vspace{0.5cm}

\centerline{\large \bf 2. Preliminaries}

\vspace{0.2cm}

Throughout the paper, $\Lambda$ is an Artinian algebra. We denote by
$(-)^*=\Hom _{\Lambda}(-, \Lambda)$ and $\mathbb{D}$ the ordinary
duality, that is, $\mathbb{D}={\rm Hom}_{R}(\ ,I^{0}(R/J(R)))$,
where $R$ is the center of $\Lambda$, $J(R)$ is the Jacobson radical
of $R$ and $I^{0}(R/J(R))$ is the injective envelope of $R/J(R)$. We
use $\mod \Lambda$ and $\gldim \Lambda$ to denote the category of
finitely generated left $\Lambda$-modules and the global dimension
of $\Lambda$, respectively. Let $M\in \mod \Lambda$. We use
$\pd_{\Lambda}M$ and $\id _{\Lambda}M$ to denote the projective and
injective dimensions of $M$, respectively. We use
$$\dots\rightarrow P_{i}(M)\rightarrow\dots\rightarrow P_{1}(M)
\rightarrow P_{0}(M)\rightarrow M\rightarrow 0$$ and $$0\rightarrow
M\rightarrow I^{0}(M)\rightarrow I^{1}(M)\rightarrow\dots\rightarrow
I^{i}(M)\rightarrow \dots $$ to denote a minimal projective
resolution and a minimal injective resolution of $M$, respectively.
Then $P_{0}(M)$ and $I^{0}(M)$ are a projective cover and an
injective envelope of $M$, respectively.

\vspace{0.2cm}

{\bf Lemma 2.1} ([IS, Proposition 1(3)]) {\it For a non-negative
integer $n$, if $\id _{\Lambda}\Lambda=\id _{\Lambda
^{op}}\Lambda=n$, then $\pd _{\Lambda}E=n$ for any non-zero direct
summand $E$ of $I^n(\Lambda)$.}

\vspace{0.2cm}

{\bf Lemma 2.2} ([M, Proposition 4.2]) {\it Assume that $\gldim
\Lambda =n (\geq 1)$. Then for any $0\neq M\in \mod \Lambda$, the
following statements are equivalent.

(1) $M$ has no non-zero projective direct summands and $\Ext
_{\Lambda}^i(M, \Lambda)=0$ for any $1 \leq i \leq n-1$.

(2) By applying the functor $(-)^*$ to a minimal projective
resolution of $M$, the induced sequence:
$$P_0(M)^* \to P_1(M)^* \to \cdots \to P_n(M)^* \to \Ext
_{\Lambda}^n(M, \Lambda) \to 0$$ is a minimal projective resolution
of $\Ext _{\Lambda}^n(M, \Lambda)$ in $\mod \Lambda ^{op}$.}

\vspace{0.2cm}

{\bf Definition 2.3} ([AuR]) A homomorphism $f: M\rightarrow N$ in
$\mod \Lambda$ is said to be {\it left minimal} if an endomorphism
$g:N\rightarrow N$ is an automorphism whenever $f=gf$. Dually, the
notion of right minimal homomorphisms is defined.

\vspace{0.2cm}

The following two properties of minimal homomorphisms are useful in
the rest of the paper.

\vspace{0.2cm}

{\bf Lemma 2.4} ([Au, Chapter II, Lemma 4.3]) {\it Let $0\rightarrow
A\stackrel{g}{\rightarrow }B\stackrel{f}{\rightarrow} C\rightarrow
0$ be a non-split exact sequence in $\mod\Lambda$.

(1) If $C$ is indecomposable, then $g:A\rightarrow B$ is left
minimal.

(2) If $A$ is indecomposable, then $f:B\rightarrow C$ is right
minimal.}

\vspace{0.2cm}

{\bf Lemma 2.5} ([HZ2, Lemma 2.5]) {\it Let $0\rightarrow
A\stackrel{g}{\rightarrow}B\stackrel{f}{\rightarrow} C\rightarrow 0$
be a non-split exact sequence in $\mod\Lambda$.

(1) If $g$ is left minimal, then $\Ext _{\Lambda}^{1}(C^{'}, A)\neq
0$ for any non-zero direct summand $C^{'}$ of $C$.

(2) If $f$ is right minimal, then $\Ext _{\Lambda}^{1}(C, A^{'})\neq
0$ for any non-zero direct summand $A^{'}$ of $A$.}

\vspace{0.2cm}

The following observation is easy to verify.

\vspace{0.2cm}

{\bf Lemma 2.6} {\it Let $f: M\rightarrow N$ be a homomorphism in
$\mod \Lambda$.

(1) If $f$ is left minimal and $N=N_1 \bigoplus N_2$ with $N_1\neq
0$, then $(1_{N_1}, 0)f\neq 0$.

(2) If $f$ is right minimal and $M=M_1 \bigoplus M_2$ with $M_1\neq
0$, then $f{\binom {1_{M_1}} 0}\neq 0$.}

\vspace{0.2cm}

{\bf Definition 2.7} ([AuR]) Let $\mathscr{C}$ be a full subcategory
of $\mod \Lambda$, $M\in\mod\Lambda$ and $C\in \mathscr{C}$. The
homomorphism $f: M\to C$ is said to be a {\it left}
$\mathscr{C}$-{\it approximation} of $M$ if ${\rm Hom}_{\Lambda}(C,
C')\to {\rm Hom}_{\Lambda}(M, C')\to 0$ is exact for any
$C'\in\mathscr{C}$. The subcategory $\mathscr{C}$ is said to be {\it
covariantly finite} in $\mod\Lambda$ if every module in
$\mod\Lambda$ has a left $\mathscr{C}$-approximation. The notions of
{\it right} $\mathscr{C}$-{\it approximations} and {\it
contravariantly finite subcategories} of $\mod\Lambda$ may be
defined dually. The subcategory $\mathscr{C}$ is said to be {\it
functorially finite} in $\mod\Lambda$ if it is both covariantly
finite and contravariantly finite in $\mod\Lambda$.

\vspace{0.2cm}

The following useful lemma is due to T. Wakamatsu.

\vspace{0.2cm}

{\bf Lemma 2.8} ([AuR, Lemma 1.3]) {\it Let $\mathscr{C}$ be a full
subcategory of $\mod \Lambda$ which is closed under extensions and
$D\in \mod \Lambda$.

(1) If $D \stackrel{f}{\rightarrow} C \rightarrow Z \rightarrow 0$
is exact with $f$ a minimal left $\mathscr{C}$-approximation of $D$,
then $\Ext _{\Lambda}^{1}(Z, \mathscr{C})=0$.

(2) If $0\rightarrow Z\rightarrow C \stackrel{f}{\rightarrow} D $ is
exact with $f$ a minimal right $\mathscr{C}$-approximation of $D$,
then $\Ext _{\Lambda}^{1}(\mathscr{C},Z)=0$.}

\vspace{0.2cm}

Let $\mathscr{C}$ be a full subcategory of $\mod \Lambda$ and $n$ a
positive integer. We denote by $^{\bot_n}\mathscr{C}= \{ X\in \mod
\Lambda \ |\ {\rm Ext}_{\Lambda}^{i}(X, C)=0$ for any $C \in
\mathscr{C}$ and $1 \leq i \leq n \}$, and $\mathscr{C}^{\bot_n}= \{
X\in \mod \Lambda \ |\ {\rm Ext}_{\Lambda}^{i}(C, X)=0$ for any $C
\in \mathscr{C}$ and $1 \leq i \leq n \}$.

\vspace{0.2cm}

{\bf Definition 2.9} ([Iy1]) Let $\mathscr{C}$ be a functorially
finite subcategory of $\mod \Lambda$. For a positive integer $n$,
$\mathscr{C}$ is called a {\it maximal} $n$-{\it orthogonal
subcategory} of $\mod\Lambda$ if
$\mathscr{C}={^{\bot_n}\mathscr{C}}=\mathscr{C}^{\bot_n}$.

\vspace{0.2cm}

From the definition above, we get easily that both $\Lambda$ and
$\mathbb{D}\Lambda ^{op}$ are in any maximal $n$-orthogonal
subcategory of $\mod \Lambda$. For a module $M \in \mod \Lambda$, we
use $\add _{\Lambda}M$ to denote the subcategory of $\mod \Lambda$
consisting of all modules isomorphic to direct summands of finite
direct sums of copies of $_{\Lambda}M$. Then
$\add_{\Lambda}(\Lambda\bigoplus \mathbb{D}\Lambda ^{op})$ is
contained in any maximal $n$-orthogonal subcategory of $\mod
\Lambda$. On the other hand, it is easy to see that if
$\add_{\Lambda}(\Lambda\bigoplus \mathbb{D}\Lambda ^{op})$ is a
maximal $n$-orthogonal subcategory of $\mod \Lambda$, then
$\add_{\Lambda}(\Lambda\bigoplus \mathbb{D}\Lambda ^{op})$ is the
unique maximal $n$-orthogonal subcategory of $\mod \Lambda$. In this
case, we say that $\Lambda$ admits a {\it trivial maximal} $n$-{\it
orthogonal subcategory} of $\mod \Lambda$.

\vspace{0.2cm}

{\bf Lemma 2.10} {\it Assume that $\gldim\Lambda=n(\geq 2)$ and
$\Lambda$ admits a maximal $(n-1)$-orthogonal subcategory of
$\mod\Lambda$. Then we have

(1) $I$ is projective or $\pd _{\Lambda}I=n$ for any injective
module $I \in \mod \Lambda$.

(2) $P$ is injective or $\id _{\Lambda}P=n$ for any projective
module $P \in \mod \Lambda$.}

\vspace{0.2cm}

{\it Proof.} Because $\Lambda$ admits a maximal $(n-1)$-orthogonal
subcategory of $\mod\Lambda$, $\Ext _{\Lambda}^i(\mathbb{D}\Lambda,
\Lambda)$\linebreak $=0$ for any $1 \leq i \leq n-1$. For an
injective module $I$ and a projective module $P$ in $\mod \Lambda$,
considering a minimal projective resolution of $I$ and a minimal
injective resolution of $P$, respectively, then it is easy to get
the assertions.  $\hfill{\square}$

\vspace{0.2cm}

Let $\mathscr{C}$ be a full subcategory of $\mod \Lambda$,
$M\in\mod\Lambda$. We use $\underline{\mathscr{C}}$ and
$\overline{\mathscr{\mathscr{C}}}$ to denote the stable
subcategories of $\mathscr{C}$ modulo projective and injective
modules, respectively. We use $\Ind \mathscr{C}$ to denote the
subcategory of $\mathscr{C}$ consisting of indecomposable modules in
$\mathscr{C}$.

\vspace{0.2cm}

{\bf Lemma 2.11} ([Iy1, Theorem 2.2.4]) {\it Assume that
$\gldim\Lambda=2$ and $\mathscr{C}$ is a maximal 1-orthogonal
subcategory of $\mod\Lambda$. Then we have

(1) There exist mutually inverse equivalences $\mathbb{D}\Ext
_{\Lambda}^{2}(-,\Lambda)$:
$\underline{\mathscr{C}}\rightarrow\overline{\mathscr{C}}$ and $\Ext
_{\Lambda}^{2}(\mathbb{D}-,\Lambda)$:
$\overline{\mathscr{C}}\rightarrow \underline{\mathscr{C}}$.

(2) $\mathbb{D}\Ext _{\Lambda}^{2}(-,\Lambda)$ gives a bijection
from the non-projective objects to the non-injective objects in
$\Ind \mathscr{C}$, and the inverse is given by $\Ext
_{\Lambda}^{2}(\mathbb{D}-,\Lambda)$.}

\vspace{0.2cm}

We recall the definition of almost hereditary algebras from [HRS].

\vspace{0.2cm}

{\bf Definition 2.12}  $\Lambda$ is called an {\it almost hereditary
algebra} if the following conditions are satisfied: (1) $\gldim
\Lambda\leq 2$; and (2) If $X \in \mod \Lambda$ is indecomposable,
then either $\id {_{\Lambda}X}\leq 1$ or $\pd {_{\Lambda}X}\leq 1$.

\vspace{0.2cm}

Also recall from [HRS] that $\Lambda$ is called a {\it quasi-tilted
algebra} if $\Lambda=\End (T)$, where $T$ is a tilting object in a
locally finite hereditary abelian $R$-category. It was proved in
[HRS, Chapter III, Theorem 2.3] that $\Lambda$ is almost hereditary
if and only if it is quasi-tilted, which implies that a tilted
algebra is almost hereditary. On the other hand, it was showed in
[HZ1, Theorem 3.15] that if an almost hereditary algebra $\Lambda$
admits a maximal 1-orthogonal subcategory $\mathscr{C}$ of
$\mod\Lambda$, then $\mathscr{C}$ is trivial (that is,
$\mathscr{C}=\add_{\Lambda}(\Lambda\bigoplus \mathbb{D}\Lambda
^{op})$).

\vspace{0.2cm}

In the subcategory of $\mod \Lambda$ consisting of indecomposable
modules, we define the relation $\rightsquigarrow _0$ given by
$X\rightsquigarrow _0Y$ if $\Hom _{\Lambda}(X, Y)\neq 0$, and let
$\rightsquigarrow$ be the transitive closure of this relation, that
is, $X\rightsquigarrow Y$ if and only if there exists a positive
integer $n$ and a path
$X=X_{0}\stackrel{f_{0}}{\rightarrow}X_{1}\stackrel{f_{1}}{\rightarrow}\dots\rightarrow
X_{n}\stackrel{f_{n}}{\rightarrow}X_{n+1}=Y$ of indecomposable
modules $X_{i}$ and non-zero homomorphisms $f_{j}:X_{j}\rightarrow
X_{j+1}$ for any $0\leq i\leq n+1$ and $0\leq j\leq n$. We use
$\mathcal {R}_{\Lambda}$ to denote the subset of indecomposable
$\Lambda$-modules given by $\mathcal {R}_{\Lambda}=\{X$ is
indecomposable in $\mod \Lambda |$ for all $X\rightsquigarrow Y$ we
have $\id_{\Lambda}Y\leq1\}$ (see [HRS]).

\vspace{0.2cm}

{\bf Lemma 2.13} {\it If $\Lambda$ is an Auslander's 1-Gorenstein
and almost hereditary algebra, then $\Lambda$ is a tilted algebra.}

\vspace{0.2cm}

{\it Proof.} By the block theory of Artinian algebras, we can assume
that $\Lambda$ is connected. Since $\Lambda$ is an almost hereditary
algebra, all indecomposable injective $\Lambda$-modules belong to
$\mathcal {R}_{\Lambda}$ by [HRS, Chapter II, Theorem 1.14] . Notice
that $\Lambda$ is an Auslander's 1-Gorenstein algebra, there exists
an indecomposable projective-injective module $I$ in $\mod \Lambda$.
So $I\in \mathcal {R}_{\Lambda}$. It follows from [HRS, Chapter II,
Corollary 3.4] that $\Lambda$ is a tilted algebra. $\hfill{\square}$

\vspace{0.5cm}

\centerline{\large \bf 3. Homological dimensions of indecomposable
modules}

\vspace{0.2cm}

In this section, we prove the main result by studying the
homological properties of simple modules and indecomposable modules
over an algebra with global dimension two and admitting a trivial
maximal 1-orthogonal subcategory.

Throughout this section, assume that $\gldim\Lambda=2$ and $\Lambda$
admits a trivial maximal 1-orthogonal subcategory
$\mathscr{C}(=\add_{\Lambda}(\Lambda\bigoplus \mathbb{D}\Lambda
^{op}))$ of $\mod\Lambda$. In this case, we then have that any
indecomposable module in $\mathscr{C}$ is either projective or
injective, and that $\Ext _{\Lambda}^1(I, P)=0$ for any injective
module $I$ and any projective module $P$ in $\mod \Lambda$. In the
following, we will use these facts freely.

The following lemma plays a crucial role in the proof of the main
result in this paper.

\vspace{0.2cm}

{\bf Lemma 3.1} {\it Let $S\in \mod \Lambda$ be a non-injective
simple module with $\pd_{\Lambda}S=2$. Then there exists a minimal
right $\mathscr{C}$-approximation of $S$: $$0\rightarrow
P\rightarrow I\rightarrow S\rightarrow 0$$ satisfying the following
conditions:

(1) $P$ is projective and $\id_{\Lambda}Q=2$ for any indecomposable
direct summand $Q$ of $P$.

(2) $I$ is injective and for any indecomposable direct summand $E$
of $I$, $\pd_{\Lambda}E=2$ and $E$ has a unique maximal submodule.}

\vspace{0.2cm}

{\it Proof.} Because $\mathscr{C}(=\add_{\Lambda}(\Lambda\bigoplus
\mathbb{D}\Lambda ^{op}))$ is a trivial maximal 1-orthogonal
subcategory of $\mod\Lambda$, $\mathscr{C}$ is functorially finite
in $\mod \Lambda$. Let $S\in \mod \Lambda$ be a non-injective simple
module with $\pd_{\Lambda}S=2$. Then there exists a minimal right
$\mathscr{C}$-approximation of $S$: $$0\rightarrow P\rightarrow
I\rightarrow S\rightarrow 0.$$

(1) By Lemma 2.8(2), $\Ext _{\Lambda}^1(\mathscr{C}, P)=0$ and $P\in
{\mathscr{C}^{\bot _1}}(=\mathscr{C})$. So, by Lemma 2.5(2), for any
indecomposable direct summand $Q$ of $P$, $Q$ is projective, but not
injective. Then it follows from Lemma 2.10(2) that
$\id_{\Lambda}Q=2$.

(2) Let $E$ be any indecomposable direct summand of $I$ such that
$I=I^{'}\bigoplus E$. Then by Lemma 2.6(2), we get the following
commutative diagram with exact columns and rows:

$$\xymatrix{ &0\ar[d]&0\ar[d]& & \\
0\ar[r]&\Ker f{\binom {1_{I^{'}}} 0}\ar[r]\ar[d]&
I^{'}\ar[r]^{f{\binom {1_{I^{'}}} 0}}\ar[d]^{\binom {1_{I^{'}}} 0}
&S\ar[r]\ar@{=}[d]&0\\
0\ar[r]&P\ar[r]\ar[d]&I\ar[r]^{f}\ar[d]&S\ar[r]&0\\
&E\ar[d]\ar@{=}[r]&E\ar[d]& & \\
&0&0& & }$$ Note that $P$ is projective by (1). We claim that $E$ is
not projective. Otherwise, if $E$ is projective, then $\Ker f{\binom
{1_{I^{'}}} 0}$ is projective and the upper row in the above diagram
is a right $\mathscr{C}$-approximation of $S$, which contradicts
with that the middle row in the above diagram is a minimal right
$\mathscr{C}$-approximation of $S$. The claim is proved. Then we
conclude that $E$ is injective and $I$ is also injective. It follows
from Lemma 2.10(1) that $\pd_{\Lambda}E=2$.

Consider the following commutative diagram with exact columns and
rows:
$$\xymatrix{ &0\ar[d]&0\ar[d]& & \\
0\ar[r]&\Ker f{\binom {1_{E}} 0}\ar[r]\ar[d]&
E\ar[r]^{f{\binom {1_{E}} 0}}\ar[d]^{\binom {1_{E}} 0}&S\ar[r]\ar@{=}[d]&0\\
0\ar[r]&P\ar[r]\ar[d]&I\ar[r]^{f}\ar[d]&S\ar[r]&0\\
&I^{'}\ar[d]\ar@{=}[r]&I^{'}\ar[d]& & \\
&0&0& & }$$ Since $E$ is indecomposable and injective, $\Ker
f{\binom {1_{E}} 0}$ is indecomposable. Because $P$ is projective by
(1) and $I^{'}$ is not projective, the leftmost column in the above
diagram is non-split and then it is a projective cover of $I^{'}$ by
Lemma 2.4(2). It follows that $P \bigoplus P_0(E)\cong P_0(I^{'})
\bigoplus P_0(E)\cong P_0(I)$. On the other hand, we have an
epimorphism $P \bigoplus P_0(S)\to I \to 0$. But $P_0(S)$ is
indecomposable, $P_0(E)\cong P_0(S)$, which implies $E$ has a unique
maximal submodule. $\hfill{\square}$

\vspace{0.2cm}

In the following, we will give some applications of Lemma 3.1, which
are needed later. We first have the following

\vspace{0.2cm}

{\bf Lemma 3.2} {\it $\id_{\Lambda}S\leq1$ for any simple module $S
\in \mod \Lambda$ with $\pd_{\Lambda}S=2$.}

\vspace{0.2cm}

{\it Proof.} If there exists a simple module $S\in \mod \Lambda$
such that $\pd_{\Lambda}S=\id_{\Lambda}S=2$, then by Lemma 3.1,
there exists a minimal right $\mathscr{C}$-approximation of $S$:
$$0\to P\buildrel {g} \over \longrightarrow I\to  S\to 0$$ with $I$ injective.
By Lemma 2.4(1), $g$ is left minimal. So
$$0\to P\buildrel {g} \over \longrightarrow I\to I^{0}(S)\to I^{1}(S)\to
I^{2}(S)\to0$$ is a minimal injective resolution of $P$ and hence
$\id _{\Lambda}P=3$, which contradicts with $\gldim \Lambda=2$.
$\hfill{\square}$

\vspace{0.2cm}

The following result is another application of Lemma 3.1.

\vspace{0.2cm}

{\bf Lemma 3.3} {\it Let $S\in \mod \Lambda$ be a simple module with
$\pd_{\Lambda}S=2$ and $\id_{\Lambda}S=1$. Then we have the
following

(1) $\Ext _{\Lambda}^{1}(S,\Lambda)\in \mod \Lambda ^{op}$ is
projective.

(2) $\Ext _{\Lambda}^{2}(S,\Lambda)\in \mod \Lambda ^{op}$ is
injective and $\pd_{\Lambda ^{op}}\Ext
_{\Lambda}^{2}(S,\Lambda)=2$.}

\vspace{0.2cm}

{\it Proof.} Let $S\in \mod \Lambda$ be a simple module with
$\pd_{\Lambda}S=2$ and $\id_{\Lambda}S=1$. By Lemma 3.1, there
exists a minimal right $\mathscr{C}$-approximation of $S:$ $$0\to
P\to I\to S\to 0$$ such that $P(\neq 0)$ is projective, $I$ is
injective without non-zero projective direct summands and $\pd
_{\Lambda}I=2$.

By Lemma 2.2, $0\to I^{*} \to P_0(I)^* \to P_1(I)^* \to P_2(I)^* \to
\Ext _{\Lambda}^2(I, \Lambda) \to 0$ is a minimal projective
resolution of $\Ext _{\Lambda}^2(I, \Lambda)$ in $\mod \Lambda
^{op}$. Because $\gldim \Lambda =2$, $I^{*}=0$. Then we have that
$\Ext _{\Lambda}^{1}(S,\Lambda)\cong P^{*}$ is projective. On the
other hand, by Lemma 2.11, we have that $\Ext
_{\Lambda}^{2}(S,\Lambda)\cong \Ext _{\Lambda}^{2}(I,\Lambda)$ is
injective and $\pd_{\Lambda ^{op}}\Ext _{\Lambda}^{2}(S,\Lambda)=2$.
$\hfill{\square}$

\vspace{0.2cm}

As a consequence of Lemma 3.3 and the results obtained in Section 2,
we have the following result.

\vspace{0.2cm}

{\bf Lemma 3.4} {\it Let $S\in \mod \Lambda$ be a simple module with
$\pd_{\Lambda}S=2$ and $\id_{\Lambda}S=1$. Then
$$0\to\Ext _{\Lambda}^{1}(S,\Lambda)\to
\Ext _{\Lambda}^{2}(I^{1}(S),\Lambda)\to\Ext
_{\Lambda}^{2}(I^{0}(S),\Lambda)\to\Ext _{\Lambda}^{2}(S,\Lambda)\to
0$$ is a minimal injective resolution of $\Ext
_{\Lambda}^{1}(S,\Lambda)$ in $\mod \Lambda ^{op}$.}

\vspace{0.2cm}

{\it Proof.} Let $S\in \mod \Lambda$ be a simple module with
$\pd_{\Lambda}S=2$ and $\id_{\Lambda}S=1$. By applying the functor
$(-)^*$ to a minimal injective resolution of $S$, we get an exact
sequence in $\mod \Lambda ^{op}$:
$$0\to\Ext _{\Lambda}^{1}(S,\Lambda)\buildrel {\alpha} \over
\longrightarrow \Ext _{\Lambda}^{2}(I^{1}(S),\Lambda)\to\Ext
_{\Lambda}^{2}(I^{0}(S),\Lambda)\to\Ext _{\Lambda}^{2}(S,\Lambda)\to
0.$$

Because $S\in \mod \Lambda$ is simple with $\pd _{\Lambda}S=2$,
$\Hom _{\Lambda}(S, I^2(\Lambda))\cong \Ext _{\Lambda}^2(S,
\Lambda)\neq 0$ and $I^0(S)$ is isomorphic to a direct summand of
$I^2(\Lambda)$. Because $\gldim \Lambda =2$, $\pd
_{\Lambda}I^0(S)=2$ by Lemma 2.1. Then $\Ext
_{\Lambda}^{2}(I^{0}(S),\Lambda)\in \mod \Lambda ^{op}$ is
indecomposable and injective by Lemma 2.11. On the other hand, $\pd
_{\Lambda}I^1(S)\neq 1$ by Lemma 2.10, so $\pd _{\Lambda}I^1(S)=2$.
In addition, by Lemma 2.5(1), $I^1(S)$ has no non-zero projective
direct summands. It follows from Lemma 2.11 that $\Ext
_{\Lambda}^{2}(I^{1}(S),\Lambda)\in \mod \Lambda ^{op}$ is injective
without non-zero projective direct summands.

By Lemma 3.3, $\Ext _{\Lambda}^{1}(S,\Lambda)$ is projective, so
$\alpha$ is non-split. Because $\Ext
_{\Lambda}^{2}(I^{0}(S),\Lambda)$ is indecomposable and injective,
$\Coker \alpha$ is indecomposable. Then by Lemma 2.4(1), $\alpha$ is
left minimal. Note that $\Ext _{\Lambda}^{2}(S,\Lambda)\in \mod
\Lambda ^{op}$ is injective by Lemma 3.2(2). The proof is finished.
$\hfill{\square}$

\vspace{0.2cm}

Now we can classify the simple modules over Auslander's 1-Gorenstein
algebras in terms of the projective and injective dimensions of
simple modules as follows.

\vspace{0.2cm}

{\bf Proposition 3.5} {\it Let $\Lambda$ be an Auslander's
1-Gorenstein algebra and $S\in \mod \Lambda$ a simple module. Then
we have

(1) $S$ is injective if $\pd_{\Lambda}S=2$.

(2) $\pd_{\Lambda}S=1$ if and only if $\id_{\Lambda}S=1$.}

\vspace{0.2cm}

{\it Proof.} (1) Let $S\in \mod \Lambda$ be a simple module with
$\pd_{\Lambda}S=2$. If $S$ is not injective, then by Lemma 3.2,
$\id_{\Lambda}S=1$. It follows from Lemmas 3.3(1) and 3.4 that
$$0\to\Ext _{\Lambda}^{1}(S,\Lambda)\to\Ext _{\Lambda}^{2}(I^{1}(S),\Lambda)
\to\Ext _{\Lambda}^{2}(I^{0}(S),\Lambda)\to\Ext
_{\Lambda}^{2}(S,\Lambda)\to 0$$ is a minimal injective resolution
of the projective module $\Ext _{\Lambda}^{1}(S,\Lambda)$ in $\mod
\Lambda ^{op}$. Because $\Lambda$ is Auslander's 1-Gorenstein,
$I^{0}(\Lambda ^{op})$ is projective, which implies that $\Ext
_{\Lambda}^{2}(I^{1}(S),\Lambda)$ is also projective. On the other
hand, by Lemma 2.10, $\pd_{\Lambda}I^{1}(S)\neq 1$, so
$\pd_{\Lambda}I^{1}(S)=2$ and hence $\pd _{\Lambda ^{op}}\Ext
_{\Lambda}^{2}(I^{1}(S),\Lambda)=2$ by Lemma 2.11, which is a
contradiction.

(2) If $\id_{\Lambda}S=1$, then $\pd_{\Lambda}S\neq 2$ by (1) and
$S$ is not projective by Lemma 2.10. So $\pd_{\Lambda}S=1$.
Conversely, if $\pd_{\Lambda}S=1$, then $\id_{\Lambda
^{op}}\mathbb{D}S=1$. Note that $\Lambda ^{op}$ is also Auslander's
1-Gorenstein and $\Lambda^{op}$ admits a trivial maximal
1-orthogonal subcategory of $\mod\Lambda^{op}$. By using an argument
dual to above, we have that $\pd_{\Lambda ^{op}}\mathbb{D}S=1$ and
$\id_{\Lambda}S=1$. $\hfill{\square}$

\vspace{0.2cm}

To prove the main result, we also need the following lemma.

\vspace{0.2cm}

{\bf Lemma 3.6} {\it Let $\Lambda$ be an Auslander's 1-Gorenstein
algebra. If both $M$ and $P$ are indecomposable modules in $\mod
\Lambda$ such that $M$ is non-projective and $P$ is projective with
$\id_{\Lambda}P=2$, then $\Hom _{\Lambda}(M, P)=0$.}

\vspace{0.2cm}

{\it Proof.} Let both $M$ and $P$ be indecomposable modules in $\mod
\Lambda$ such that $M$ is non-projective and $P$ is projective with
$\id_{\Lambda}P=2$. If there exists a non-zero homomorphism
$f\in\Hom _{\Lambda}(M,P)$, then $f$ is not epic. Since $P$ is
indecomposable and projective, $J(\Lambda)P$ is the unique maximal
submodule of $P$ (where $J(\Lambda)$ is the Jacobson radical of
$\Lambda$) and $P/J(\Lambda)P(=S)$ is a simple module. Consider the
following diagram:
$$\xymatrix{\ & &M\ar[d]^{f}{\ar@{-->}[ld]_{f^{'}}}& & \\
0\ar[r]&J(\Lambda)P\ar[r]^{i}&P\ar[r]^{\pi}&S\ar[r]&0 }$$ where $i$
is the inclusion homomorphism and $\pi$ is the natural epimorphism
(which is a projective cover of $S$). Because $f(\neq 0)$ is not
epic, there exists a non-zero homomorphism $f^{'}$ such that
$f=if^{'}$.

If $S$ is projective, then $P(\cong S)$ is simple and so $f$ is
epic, which is a contradiction. If $\pd _{\Lambda}S=2$, then by
Proposition 3.5, $S$ is injective and $\mathbb{D}S$ is projective.
On the other hand, $\mathbb{D}\pi: \mathbb{D}S\hookrightarrow
\mathbb{D}P$ is an injective envelope of $\mathbb{D}S$ and
$\pd_{\Lambda^{op}}\mathbb{D}P=2$, which is also a contradiction
because $\Lambda$ (and hence $\Lambda ^{op}$) is Auslander's
1-Gorenstein. Thus we conclude that $\pd_{\Lambda}S=1$. So
$J(\Lambda)P$ is projective and hence the injective dimension of any
indecomposable direct summand of $J(\Lambda)P$ is equal to 2 by
Lemmas 2.5(2) and 2.10. Because $f^{'}\neq 0$, there exists a
non-zero indecomposable direct summand $P_{1}$ of $J(\Lambda)P$ such
that $\Hom _{\Lambda}(M, P_1)\neq 0$. Then iterating the above
process, we get an infinite descending chain of non-zero submodules
of $P$:
$$P\supsetneqq P_{1}\supsetneqq P_{2}\supsetneqq\cdots,$$
which is a contradiction. The proof is finished. $\hfill{\square}$

\vspace{0.2cm}

Now we are in a position to state the main result in this paper.

\vspace{0.2cm}

{\bf Theorem 3.7} {\it Let $\Lambda$ be an Auslander's 1-Gorenstein
algebra. Then for any indecomposable module $M\in\mod\Lambda$, we
have

(1) $\pd_{\Lambda}M=1$ if and only if  $\id_{\Lambda}M=1$.

(2) $M$ is injective if $\pd_{\Lambda}M=2$.}

\vspace{0.2cm}

{\it Proof.} Let $M\in\mod\Lambda$ be an indecomposable and
non-projective module. We claim that $\id_{\Lambda}M\leq 1$. Take a
minimal left $\mathscr{C}$-approximation of $M$:
$$0\to M\to C_{0}\to C_{1}\to 0.$$ If $P$ is a non-zero
indecomposable and projective direct summand of $C_0$, then by Lemma
2.6(1), $\Hom _{\Lambda}(M, P)\neq 0$. So $\id _{\Lambda}P\neq 2$ by
Lemma 3.6 and hence $P$ is injective by Lemma 2.10. It follows that
$C_0$ is injective.

If $C_1=0$, then $M\cong C_0$ is injective (in this case,
$\pd_{\Lambda}M=2$ by Lemma 2.10). Now suppose $C_1 \neq 0$. By
Lemma 2.8(1), $C_1 \in {^{\bot _1}\mathscr{C}}(=\mathscr{C})$. On
the other hand, $C_1$ has no non-zero projective direct summands by
Lemma 2.5(1), so $C_1$ is injective and hence $\id_{\Lambda}M=1$.
The claim is proved.

(1) If $\pd_{\Lambda}M=1$, then $\id_{\Lambda}M=1$ by the above
argument and Lemma 2.10. Because $\Lambda ^{op}$ is also Auslander's
1-Gorenstein and $\Lambda^{op}$ admits a trivial maximal
1-orthogonal subcategory of $\mod\Lambda^{op}$, we get dually that
$\id_{\Lambda}M=1$ implies $\pd_{\Lambda}M=1$.

(2) By the above argument, it is easy to see that $M$ is injective
if $\pd_{\Lambda}M=2$. $\hfill{\square}$

\vspace{0.2cm}

As an immediate consequence of Theorem 3.7, we get the following
result.

\vspace{0.2cm}

{\bf Theorem 3.8} {\it If $\Lambda$ is an Auslander's 1-Gorenstein
algebra, then $\Lambda$ is a tilted algebra.}

\vspace{0.2cm}

{\it Proof.} If $\Lambda$ is an Auslander's 1-Gorenstein algebra,
then $\Lambda$ is almost hereditary by Theorem 3.7. It follows from
Lemma 2.13 that $\Lambda$ is tilted. $\hfill{\square}$

\vspace{0.2cm}

At the end of this paper, we give the following example to
illustrate the obtained results above.

\vspace{0.2cm}

{\bf Example 3.9} Let $n\geq 4$ and $\Lambda$ be a
finite-dimensional algebra given by the quiver:
$$\xymatrix{1 &
\ar[l]_{\alpha _{1}} 2 & \ar[l]_{\alpha _{2}} 3 & \ar[l]_{\alpha
_{3}} \cdots & \ar[l]_{\alpha _{n-1}} n}
$$ modulo the ideal
generated by $\{\alpha_{1}\alpha_{2}\cdots \alpha_{n-1}\}$. Then

(1) $\Lambda$ is a tilted algebra of finite representation type, and
$\Lambda$ admits a trivial maximal 1-orthogonal subcategory $\add
_{\Lambda} [(\bigoplus_{i=1}^nP(i))\bigoplus
(\bigoplus_{i=3}^nI(i))]$ of $\mod\Lambda$. Because $I^0(\Lambda)$
is projective and $\pd _{\Lambda}I^1(\Lambda)=2$, $\Lambda$ is
Auslander's 1-Gorenstein, but not an Auslander algebra.

(2) For any indecomposable module $M\in\mod\Lambda$, we have that
$\pd_{\Lambda}M=1$ if and only if  $\id_{\Lambda}M=1$ and that $M$
is injective if $\pd_{\Lambda}M=2$.

\vspace{0.5cm}

{\bf Acknowledgements} The research was partially supported by the
Specialized Research Fund for the Doctoral Program of Higher
Education (Grant No. 20060284002), NSFC (Grant No. 10771095) and NSF
of Jiangsu Province of China (Grant No. BK2007517).

\end{document}